\documentclass[preprint, times, 3p, 11pt]{elsarticle}
\usepackage{amsmath}
\usepackage{enumerate}
\usepackage{multirow}
\usepackage{graphicx}
\usepackage{amssymb}
\usepackage{amsthm}
\usepackage{epstopdf}
\usepackage{xcolor}
\usepackage{mathrsfs}
\usepackage[colorlinks=true]{hyperref}
\usepackage{doi}
\usepackage{hyperref}
\usepackage{marvosym}
\hypersetup{citecolor=red}

\usepackage{amssymb,amsthm,amsmath}

\usepackage{url}

 \newtheorem{thm}{Theorem}[section]
\newtheorem{lem}{Lemma}[section]

\theoremstyle{definition}
\newtheorem{defn}{Definition}[section]

\numberwithin{equation}{section}

\usepackage{amssymb,amsmath,color}
\numberwithin{equation}{section}

\usepackage{lineno}

\begin{document}\bibliographystyle{plain}
\begin{frontmatter}
\title{\textbf{The Shadowing Properties Of  Nonautonomous   Dynamical   System}\tnoteref{label1}}
\tnotetext[label1]{ Data: March 21, 2024}

\author[a]{Min An$^{\href{mailto:202006021023@stu.cqu.edu.cn}{\textrm{\Letter}}}$}


\address[a]{College of Mathematics and Statistics,
 Chongqing University, Chongqing 401331, China}
\begin{abstract}
Let $\left(X_n, d_n\right)$ be a sequence of metric spaces and let $\mathcal{F}=\left\{f_n\right\}_{n \in \mathbb{Z}}$ be a sequence of continuous and onto maps $f_n: X_n \rightarrow X_{n+1}, n \in \mathbb{Z}_{+}$. In this paper, we prove that if the compression ratio meets $\prod \lambda_i=0$, then there exists $\delta_n>0$ such that any $\delta_n$ - pseudo-orbit $\left\{x_n\right\}_{n \in \mathbb{Z}}$ is $\varepsilon$ - shadowed by a unique point $x \in X_0$. For the asymptotic average shadowing property, we prove that if $\left.\mathcal{F}\right|_A$ has asymptotically average shadowing property, then $\mathcal{F}$ also has a asymptotically average shadowing propertywhen a density-related condition is satisfied. Additionally, the conclusion that the shadowing performance of strong equicontinuity and pseudo-shadowing property implies limit shadowing is also obtained. Furthermore, the shadowing property of non-autonomous product space is also discussed.
\end{abstract}

\end{frontmatter}


\section{Introduction}\label{Section 1}

Pseudo-orbit shadowing is not only a very important concept in dynamical system theory, but also the most important dynamic property and technical tool in dynamical system theory. The pseudo-shadowing  property studies whether the pseudo-orbit under a mapping can be shadowed by the real orbit . Based on the needs of theory and application, different concepts of pseudo shadowing and shadowing property have been proposed from different angles and standards, such as average pseudo shadowing and asymptotic average pseudo shadowing, asymptotic average shadowing, Lipschitz shadowing, limit shadowing.

In 2019, Castro A.\cite{1pp} obtained some conclusions about expanding maps in non-autonomous dynamic systems; In 2004, Junling Zhao\cite{2pp} proved the limit shadowing property of expanding maps with equicontinuity in autonomous dynamic systems. In 2011, Kulczycki M.\cite{3pp} proved the conclusion that the asymptotic average property of $f$ can be obtained from the asymptotic average property of $\left.f\right|_A$ in autonomous dynamical systems. Additionally, many conclusions about product space in non-autonomous dynamical systems are obtained in \cite{5pp,6pp,4pp}.

The structure of this paper is as follows: In section 2, preparatory knowledge is introduced; In Section 3, we prove first that if the compression ratio meets $\prod \lambda_i=0$, then there exists $\delta_n>0$ such that any $\delta_n-$ pseudo-orbit $\left\{x_n\right\}_{n \in \mathbb{Z}_{+}}$is $\varepsilon$ - shadowed by a unique point $x \in X_0$. Secondly, we prove some conclusions about limit shadowing and asymptotic average shadowing in non-autonomous dynamical systems.; In Section 4, we discuss the shadowing peoperty of product spaces in non-autonomous dynamical systems.

\section{Preliminary}\label{Section 2}

Throughout this article, we have agreed that $\mathbb{Z}$ is all integers, $\mathbb{Z}_{+}$is the set of positive integers, and $\mathbb{N}$ is the set of natural numbers. The following definitions and theorems can be found in reference \cite{1pp,8pp,6pp}.

\begin{defn} 
Let $\left(X_n, d_n\right)$ be a sequence of metric spaces(here $d_n$ denotes the metric in $X_n$ )and let $\mathcal{F}=\left\{f_n\right\}_{n \in \mathbb{Z}_{+}}$be a sequence of continuous and onto maps $f_n: X_n \rightarrow X_{n+1}$, $n \in \mathbb{Z}_{+}, F_n:=f_{n-1} \circ \cdots \circ f_2 \circ f_1 \circ f_0, n \geq 0 . \mathcal{F}$ is called a family of time-varying mappins, $(X, d, \mathcal{F})$ is non-autonomous dynamical systems.
\end{defn}

\begin{defn} 
Let $\left(X_n, d_n\right)$ be a sequence of metric spaces and let $\mathcal{F}=\left\{f_n\right\}_{n \in \mathbb{Z}_{+}}$be a sequence of continuous and onto maps $f_n: X_n \rightarrow X_{n+1}, n \in \mathbb{Z}_{+}$. We say  $\mathcal{F}$ is a sequence of expanding maps if there are $\delta_0>0$ and a sequence $\left\{\lambda_n\right\}_{n \in \mathbb{Z}_{+}}$of constants in $(0,1)$ so that the following holds: for any $n \in \mathbb{N}, x \in X_{n+1}$ and $x_i \in f_n^{-1}(x)$, there exists a well-defined inverse branch $f_{n, x_i}^{-1}: B\left(x, \delta_0\right) \rightarrow V_{x_i}$ (open neighborhood of $\left.x_i\right)$ so that for every $y, z \in B\left(x, \delta_0\right)$,

$$
d_n\left(f_{n, x_i}^{-1}(y), f_{n, x_i}^{-1}(z)\right) \leq \lambda_n d_{n+1}(y, z).
$$

We say that $\mathcal{F}$ is a expansive sequence if there exists $\varepsilon_0>0$ (called an expansiveness constant) such that if $x, y \in X_0$ and $d_{n+1}\left(F_n(x), F_n(y)\right) \leq \varepsilon_0$ for every $N \in \mathbb{Z}_{+}$, then $x=y$.

We represent all metrics $d_n$ simply by $d$ for notational simplicity.

\end{defn}

\begin{defn} 
 Let $\left(X_n, d_n\right)$ be a sequence of metric spaces and let $\mathcal{F}=\left\{f_n\right\}_{n \in \mathbb{Z}_{+}}$be a sequence of continuous and onto maps $f_n: X_n \rightarrow X_{n+1}, n \in \mathbb{Z}_{+}$. Given $\delta>0$, we say  $\left\{x_n\right\}_{n=0}^{\infty}$ is a $\delta$-pseudo-orbit for $\mathcal{F}=\left\{f_n\right\}_{n \in \mathbb{Z}_{+}}$if $x_i \in X_i$ and $d\left(f_i\left(x_i\right), x_{i+1}\right)<\delta$ for every $i \geq 0$.
\end{defn}

\begin{defn} 
Let $\left(X_n, d_n\right)$ be a sequence of metric spaces and let $\mathcal{F}=\left\{f_n\right\}_{n \in \mathbb{Z}_{+}}$be a sequence of continuous and onto maps $f_n: X_n \rightarrow X_{n+1}, n \in \mathbb{Z}_{+}$. We say the sequence of maps $\mathcal{F}$ has the shadowing property if for every $\varepsilon>0$ there exists $\delta>0$ such that for any $\delta$-pseudo-orbit $\left\{x_n\right\}_{n=0}^{\infty}$ there exists $x \in X_0$ so that its $\mathcal{F}$-orbit $\varepsilon$-shadows the sequence $\left\{x_n\right\}_{n=0}^{\infty}$, that is, if $d\left(f_n\left(x_n\right), x_{n+1}\right)<\delta$ for $n \geq 0$, there exists $x \in X_0$ so that $d\left(F_n(x), x_n\right)<\varepsilon$ for $n \geq 0$.

\end{defn}

\begin{defn}
Let $\left(X_n, d_n\right)$ be a sequence of metric spaces and let $\mathcal{F}=\left\{f_n\right\}_{n \in \mathbb{Z}_{+}}$be a sequence of continuous and onto maps $f_n: X_n \rightarrow X_{n+1}, n \in \mathbb{Z}_{+}$. We say the sequence of maps $\mathcal{F}$ has the periodic shadowing property if for every $\varepsilon>0$ there exists $\delta>0$ such that for any $\delta$-pseudo-orbit $\left\{x_n\right\}_{n=0}^{\infty}$ satisfying $x_0=x_k$ is $\varepsilon$-shadowed by a fixed point $x \in X_0$.
\end{defn}

\begin{defn}
 Let $\left(X_n, d_n\right)$ be a sequence of metric spaces and let $\mathcal{F}=\left\{f_n\right\}_{n \in \mathbb{Z}_{+}}$be a sequence of continuous and onto maps $f_n: X_n \rightarrow X_{n+1}, n \in \mathbb{Z}_{+}$. We say the sequence of maps $\mathcal{F}$ has the Lipschitz shadowing property if there exists a uniform constant $L>0$ and $\delta_0>0, \forall 0<\delta<\delta_0$ so that any $\delta$-pseudo-orbit $\left\{x_n\right\}_{n=0}^{\infty}$ of $\mathcal{F}$ is $L \varepsilon$-shadowed by a point, $x \in X_0$.
\end{defn}

\begin{thm}\label{thm:2.1}
 (1) Let $\left(X_n, d_n\right)$ be a sequence of metric spaces and let $\mathcal{F}=\left\{f_n\right\}_{n \in \mathbb{Z}_{+}}$ be a sequence of expanding maps and onto maps $f_n: X_n \rightarrow X_{n+1}, n \in \mathbb{Z}_{+}$. If the contraction rates $\left\{\lambda_n\right\}_{n \in \mathbb{Z}_{+}}$satisfy $\sup _{n \in \mathbb{Z}_{+}} \lambda_n<1$ then the following hold:

(i) $\mathcal{F}$ satisfies the Lipschitz shadowing property.

(ii) For any for every $\varepsilon>0$ there exists $\delta>0$ so that any $\delta$-pseudo-orbit $\left\{x_n\right\}_{n=0}^{\infty}$ is $\varepsilon$-shadowed by a unique point $x \in X_0$.

(iii) If, in addition, $X_n=X_0$ for all $n \in \mathbb{Z}_{+}$then $\mathcal{F}$ satisfies the periodic shadowing property.
\end{thm}

\begin{defn}
 
 Let $\left(X_n, d_n\right)$ be a sequence of metric spaces and let $\mathcal{F}=\left\{f_n\right\}_{n \in \mathbb{Z}_{+}}$be a sequence of continuous and onto maps $f_n: X_n \rightarrow X_{n+1}, n \in \mathbb{Z}_{+}$. Let $f_i^n=f_n \circ f_{n-1} \circ \ldots \circ f_i$. We say the sequence of maps $\mathcal{F}$ is strong-equicontinuous if for every $\varepsilon>0$ there exists $\delta>0$ so that $d(x, y)<\delta$, then
 
 $$
d\left(f_m^n(x), f_m^n(y)\right)<\varepsilon, n \geq m \geq 1.
$$

\end{defn}

\begin{defn}

Let $\left(X_n, d_n\right)$ be a sequence of metric spaces and let $\mathcal{F}=\left\{f_n\right\}_{n \in \mathbb{Z}_{+}}$be a sequence of continuous and onto maps $f_n: X_n \rightarrow X_{n+1}, n \in \mathbb{Z}_{+}$. We say the sequence of maps $\mathcal{F}$ has the limit shadowing property if for any sequence $\left\{x_n\right\}_{n \in \mathbb{Z}_{+}}$satisfy
$$
\lim _{i \rightarrow+\infty} d\left(f_i\left(x_i\right), x_{i+1}\right)=0.
$$

\noindent there exists $y \in X_0$ so that $\lim \limits_{i \rightarrow+\infty} d\left(F_i(y), x_i\right)=0$.
\end{defn}

\begin{thm}\label{thm:2.2}
(4) Let $(X, d)$ be metric space and $f$ is asurjection on $X$. If equicontinuous map $f$ has shadowing property, then $f$ has the limit shadowing property.
\end{thm}

\begin{defn}
 Let $\left(X_n, d_n\right)$ be a sequence of metric spaces and let $\mathcal{F}=\left\{f_n\right\}_{n \in \mathbb{Z}_{+}}$ be a sequence of continuous and onto maps $f_n: X_n \rightarrow X_{n+1}, n \in \mathbb{Z}_{+}$. A sequence $\left\{x_n\right\}_{n \in \mathbb{Z}_{+}}$is called an asymptotic-average pseudo orbit of $\mathcal{F}$ if
$$
\lim _{i \rightarrow+\infty} \frac{1}{n} \sum_{i=0}^{n-1} d\left(f_i\left(x_i\right), x_{i+1}\right)=0 .
$$

A sequence $\left\{x_n\right\}_{n \in \mathbb{Z}_{+}}$is said to be asymptotically shadowed in average by the point $z$ in $X_0$ if
$$
\lim _{i \rightarrow+\infty} \frac{1}{n} \sum_{i=0}^{n-1} d\left(F_i(z), x_i\right)=0 .
$$

We say $\mathcal{F}$ has asymptotic average shadowing property if any asymptotic-average pseudo orbit of $\mathcal{F}$, asymptotically shadowed in average by some point $z$ in $X_0$.

\end{defn}

\begin{thm}\label{thm:2.3}

 (5)Let $(X, d)$ be metric space and $f$ is homeomorphism on $X . A$ is invariant closed subset in $X$. If $\left.f\right|_A$ has asymptotic average shadowing property and for any $\varepsilon>0$ there exists $n \in \mathbb{N}$ so that for every $x \in X$,
$$
\frac{1}{n} \#\left\{0 \leq i<n: d\left(f^i(x), A\right)<\varepsilon\right\}>1-\varepsilon.
$$

\noindent  Then $f$ has asymptotic average shadowing property.

\end{thm}

\section{The Shadowing Properties of Nonautonomous Dynamical System}\label{Section 2}

Based on several theorems in Chapter 2, we prove conclusions about shadowing in non-autonomous dynamical systems in sections 3.1 and 3.2.

\subsection{Description of Theorem}

\begin{thm}\label{thm:3.1}
 Let $\left(X_n, d_n\right)$ be a sequence of metric spaces and let $\mathcal{F}=\left\{f_n\right\}_{n \in \mathbb{Z}_{+}}$be a sequence of expanding maps and onto maps $f_n: X_n \rightarrow X_{n+1}, n \in \mathbb{Z}_{+}$. If the contraction rates $\left\{\lambda_n\right\}_{n \in \mathbb{Z}_{+}}$satisfy $\prod\limits_{i=0}^{\infty} \lambda_i=0$, then the following hold:
 
(i) For any for every $\varepsilon>0$ there exists $\delta_n>0$ so that any $\delta_n$-pseudo-orbit $\left\{x_n\right\}_{n=0}^{\infty}$ is $\varepsilon$-shadowed by a unique point $x \in X_0$.

(ii) If $X_n=X_0$ for all $n \in \mathbb{Z}_{+}$, then $\mathcal{F}$ satisfies the periodic shadowing property.
\end{thm}

\begin{thm}\label{thm:3.2}
 Let $\left(X_n, d_n\right)$ be a sequence of compact metric spaces and let $\mathcal{F}=\left\{f_n\right\}_{n \in \mathbb{Z}_{+}}$be a sequence of continuous and onto maps $f_n: X_n \rightarrow X_{n+1}, n \in \mathbb{Z}_{+} \cdot \mathcal{F}$ satisfy: there are $\delta_0>0$ so that for any $n \in \mathbb{N}, w \in X_{n+1}$ and $z \in f_n^{-1}(w)$, there exists a well-defined inverse branch $f_{n x_i}^{-1}: B\left(w, \delta_0\right) \rightarrow V_z$ (open neighborhood of $x_i$ ). If equicontinuous $\mathcal{F}$ has shadowing property, then $\mathcal{F}$ has limit shadowing property.
\end{thm}

\textbf{Proof of Theorem 3.1.}

 (i) Let $\mathcal{F}=\left\{f_n\right\}_{n \in \mathbb{Z}_{+}}$be a sequence of expanding maps as above and let $\delta_0>0$ be so that for any $n \in \mathbb{Z}_{+}, w \in X_{n+1}$ and $z \in f_n^{-1}(w)$,there exists a well defined inverse branch $f_{n, x_i}^{-1}: B\left(w, \delta_0\right) \rightarrow V_z$ (open neighborhood of $\left.x_i\right)$. By the difinition of $\lambda_n$, for any $y_1, y_2 \in B\left(w, \delta_0\right)$, the following is established:
$$
d_n\left(f_{n, z}^{-1}\left(y_1\right), f_{n, z}^{-1}\left(y_2\right)\right) \leq \lambda_n d\left(y_1, y_2\right).
$$

Take $0<\varepsilon<\delta_0 / 2$ arbitrary and $\delta_n$ satisfy
$$
0<\delta_n<\left(1-\lambda_n\right) \varepsilon .
$$

\noindent  Given $k \geq 1$ and an arbitrary finite $\delta_n$-pseudo $\left\{x_n\right\}_{n=0}^k$, we proceed to prove that there exists $x \in X_0$ that $\varepsilon$-shadows $\left\{x_n\right\}_{n=0}^k$.

\noindent Since $0<\varepsilon<\delta_0 / 2$, inverse branches are well defined in balls of radius $\varepsilon .\left\{f_n\right\}_{n \in \mathbb{Z}_{+}}$is expanding, so for $n \in \mathbb{Z}_{+}, w \in X_{n+1}$ and $z \in f_n^{-1}(w)$, we have
$$
d\left(z, f_{n, z}^{-1}(y)\right) \leq \lambda_n d(w, y) \leq \lambda_n \varepsilon,
$$
so $f_{n, z}^{-1}(B(w, \varepsilon)) \subset B\left(z, \lambda_n \varepsilon\right)$. Next, take $\delta_{\mathrm{n}}>0$ such that
$$
\delta_n+\lambda_n \varepsilon<\varepsilon<\delta_0 / 2, \forall n \in \mathbb{Z}_{+} \cdot
$$

By $d\left(f_n\left(x_n\right), x_{n+1}\right)<\delta_n, \forall n \in \mathbb{N}$ and the triangular, we get

(i) $d\left(f_{n, x_n}^{-1}\left(B\left(f_n\left(x_n\right), \varepsilon\right)\right)\right) \leq \lambda_n d\left(B\left(f_n\left(x_n\right), \varepsilon\right)\right), \quad 0 \leq n \leq k-1$.

(ii) $f_{n, x_n}^{-1}\left(B\left(f_n\left(x_n\right), \varepsilon\right)\right) \subset B\left(f_{n-1}\left(x_{n-1}\right), \delta+\lambda_n \varepsilon\right) \subset B\left(f_{n-1}\left(x_{n-1}\right), \varepsilon\right), \quad 0 \leq n \leq k-1$.

Let
$$
Y_{k-1}^{(k)}:=f_{k-1, x_{k-1}}^{-1}\left(\overline{B\left(x_k, \varepsilon\right)}\right) \cap \overline{B\left(f_{k-1}\left(x_{k-1}\right), \varepsilon\right)},
$$

\noindent and

$$
Y_{k-j-1}^{(k)}:=f_{k-j-1, x_{k-j-1}}^{-1}\left(\overline{B\left(f\left(x_{k-j-1}\right), \varepsilon\right.}\right) \cap Y_{k-j}^{(k)}, \quad 1 \leq j \leq k-1 .
$$

\noindent By definition, $Y_0^{(k)}$ is a non-empty closed set of diameter at most $2 \varepsilon \prod_{i=1}^k \lambda_i$ and
$$
d\left(F_n(x), x_n\right)<\varepsilon, \quad 0 \leq n \leq k, \forall x \in Y_0^{(k)} .
$$

\noindent This proves arbitrary finite $\delta_n$-pseudo $\left\{x_n\right\}_{n=0}^k$ is shadowed by a point $x \in Y_0^{(k)} \subset X_0$.

$k \geq 1$ in the above discussion is arbitrarily given, so for any $\varepsilon>0$, we can get a sequence $\left\{\delta_n\right\}_{n=1}^{+\infty}$ such that for any infinite $\delta_n$-pseudo $\left\{x_n\right\}_{n=0}^{+\infty}$, there exists a sequence closed set $\left\{Y_0^{(k)}\right\}_{k=1}^{+\infty}$ in $X_0$ satisfy $Y_0^{(k+1)} \subset Y_0^{(k)}, \forall k \geq 1$. On the other hand, we observe that $\operatorname{diam}\left(Y_0^{(k)}\right)$ satisfies.
$$
\operatorname{diam}\left(Y_0^{(k)}\right) \leq 2 \varepsilon \prod_{i=1}^k \lambda_i \rightarrow 0, k \rightarrow+\infty .
$$

\noindent By compactness of $X_0, \bigcap_{k \geq 1} Y_0^{(k)}$ is singleton set. Let $\{x\}=\bigcap_{k \geq 1} Y_0^{(k)}$. Obviously, $x$ is the unique point of shadowing.

(ii) Let $\left\{x_i\right\}_{i=0}^{+\infty}$ is $n$-periodic $\delta_i$-pseudo of $\mathcal{F}$, that is
$$
d\left(f_{i-1}\left(x_{i-1}\right), x_i\right) \leq \delta_i, \quad i \in \mathbb{Z}_{+}
$$

\noindent and

$$
x_{k n+j}=x_j, \quad \forall k>0,0 \leq j<n .
$$
$\left\{x_i\right\}_{i=0}^{+\infty}$ is $\delta_i$-pseudo of $\mathcal{F}$, on the basis of (i), there exists unique point $x \in X_0 \varepsilon$-shadowing $\left\{x_i\right\}_{i=0}^{+\infty}$. That is
$$
d\left(F_{i-1}(x), x_i\right) \leq \varepsilon, \quad i \in \mathbb{Z}_{+} .
$$

Since $\left\{x_i\right\}_{i=0}^{+\infty}$ is $n$-periodic, there are
$$
x_0, x_1, \cdots, x_{n-1}, x_n, \cdots
$$

\noindent and

$$
x_n, x_{n+1}, \cdots, x_{2 n-1}, x_{2 n}, \cdots.
$$

\noindent They are the same sequence, so
$$
d\left(F_{n+i-1}(x), x_i\right) \leq \varepsilon, i \in \mathbb{Z}_{+} .
$$

This proves that $F_n(x)$ is also the $\varepsilon$-shadowing pseudo $\left\{x_i\right\}_{i=0}^{+\infty}$. By uniqueness of shadowing property, we get $F_n(x)=x$. Hence, $\mathcal{F}$ satisfies the periodic shadowing property.

\textbf{Proof of  Theorem 3.2.}

 By shadowing property of $\mathcal{F}$, we get for any $n \in \mathbb{Z}_{+}$, there exists $\delta_n>0$ such that for any $\delta_n$-pseudo $\left\{x_i\right\}_{i \in \mathbb{Z}_{+}}$of $\mathcal{F}$, there is $y \in X_0$ satisfying
$$
d\left(F_i(y), x_i\right)<\frac{1}{n}, \quad \forall i \in \mathbb{Z}_{+} .
$$

Let $\left\{x_n\right\}_{n \in \mathbb{Z}_{+}}$satisfy
$$
d\left(f_i\left(x_i\right), x_{i+1}\right) \rightarrow 0,(i \rightarrow+\infty).
$$

\noindent We next prove there exists $y \in X_0$ such that

$$
d\left(F_i(y), x_i\right) \rightarrow 0,(i \rightarrow+\infty) .
$$

In fact, according $d\left(f_i\left(x_i\right), x_{i+1}\right) \rightarrow 0(i \rightarrow+\infty)$, we get that for $n \in \mathbb{Z}_{+}$, there exists $k_n \in \mathbb{Z}_{+}$such that

$$
d\left(f_i\left(x_i\right), x_{i+1}\right)<\delta_n, \forall i \geq k_n .
$$

\noindent For every $i=0,1,2, \cdots k_n-1$, we can find $x_{n i} \in X_i$ such that
$$
f_i\left(x_{n i}\right)= \begin{cases}x_{n(i+1)}, & i=0,1,2, \cdots k_n-2, \\ x_i, & i=k_n-1 .\end{cases}
$$

\noindent Hence,
$$
x_{n 0}, x_{n 1}, \cdots x_{n\left(k_n-1\right)}, x_{k_n}, x_{k_n+1}, \cdots
$$

\noindent is $\delta_n$-pseudo of $\mathcal{F}$. Obviously, for $i=0,1,2, \cdots k_n-1$, the following is established:
$$
d\left(F_i\left(y_n\right), x_i\right)=d\left(f_{i-1} \circ \cdots \circ f_0\left(y_n\right), x_i\right)=d\left(x_i, x_i\right)=0 .
$$

\noindent Additionally, by the shadowing property of $\mathcal{F}$, there exists $y_n \in X_0$, for $i \geq k_n$, we have
$$
d\left(F_i\left(y_n\right), x_i\right)<\frac{1}{n} .
$$

\noindent That is

$$
d\left(F_i\left(y_n\right), x_i\right)<\frac{1}{n}, \quad \forall i \in \mathbb{Z}_{+} .
$$

Since $X$ is compact, let $y_n \rightarrow y(n \rightarrow+\infty)$ (otherwise, take a convergent subsequence of. Similar is discussed below). For any $\varepsilon>0, \mathcal{F}$ is strong-equicontinuous, so there is $\delta>0$ such that $\forall i \in \mathbb{Z}_{+}, \forall u, v \in X_0$, if $d(u, v)<\delta$, then
$$
d\left(f_{i-1} \circ \ldots \circ f_0(u), f_{i-1} \circ \ldots \circ f_0(v)\right)<\frac{\varepsilon}{2},
$$

\noindent that is

$$
d\left(F_i(u), F_i(v)\right)<\frac{\varepsilon}{2} .
$$

\noindent Because $y_n \rightarrow y(n \rightarrow+\infty)$, take $n \in \mathbb{Z}_{+}$make $\frac{1}{n}<\frac{\varepsilon}{2}$ and $d\left(y_n, y\right)<\delta$. By the strong-equicontinuity of $\mathcal{F}$, we get
$$
d\left(F_i\left(y_n\right), F_i(y)\right)<\frac{\varepsilon}{2}, \quad \forall i \geq 1 .
$$

\noindent On the other hand, when $i \geq k_n$,
$$
d\left(F_i\left(y_n\right), x_i\right)<\frac{1}{n}<\frac{\varepsilon}{2}.
$$

\noindent So when $i \geq k_n$, we have
$$
\begin{aligned}
d\left(F_i(y), x_i\right) & \leq d\left(F_i(y), F_i\left(y_n\right)\right)+d\left(F_i\left(y_n\right), x_i\right) \\
& <\frac{\varepsilon}{2}+\frac{\varepsilon}{2} \\
& =\varepsilon .
\end{aligned}
$$

\noindent This means
$$
d\left(F_i(y), x_i\right) \rightarrow 0,(i \rightarrow+\infty) .
$$

\noindent Hence, $\mathcal{F}$ has limit shadowing property.

\subsection{Asymptotic Average Shadowing}

\begin{thm}\label{thm:3.3}
 
 Let $(X, d)$ be compact metric space and $\mathcal{F}=\left\{f_n\right\}_{n \in \mathbb{Z}_{+}}$be a sequence of continuous and onto maps $f_n: X \rightarrow X, n \in \mathbb{Z}_{+}. A$ is invariant closed subset in $X$. If $\left.\mathcal{F}\right|_A$ has asymptotic average shadowing property and for any $\varepsilon>0$ there exists $n \in \mathbb{N}$ so that for every $x \in X$,
$$
\frac{1}{n} \#\left\{0 \leq i<n: d\left(F_i(x), A\right)<\varepsilon\right\}>1-\varepsilon,
$$

\noindent then $\mathcal{F}$ has asymptotic average shadowing property.
 
\end{thm}

\begin{lem}\label{lem:3.1}

 (5)Let $J_i \subset \mathbb{N}$ be a sequence such that
$$
\lim _{i \rightarrow+\infty} \lim _{k \rightarrow+\infty} \sup \frac{1}{k} \#\left(J_i \cap[0, k]\right)=0 .
$$

\noindent Then there is a set $J$ of density zero and increasing sequences $\left\{m_i\right\}_{i=0}^{\infty}$ and $\left\{l_i\right\}_{i=0}^{\infty}$ such that $m_0=0$ and for every $i=1,2 \ldots$ we have
$$
J \cap\left[m_{i-1}, m_i\right)=J_{l_i} \cap\left[m_{i-1}, m_i\right).
$$

\noindent Additionally, the following conditions hold:

(i) when each set $J_i$ is of density zero we can take $l_i=i$ for all $i$;

(ii) for any given sequence of infinite sets $R_i \subset \mathbb{N}$ we can choose $m_i$ in such a way that
$
m_i \in R_i.
$

\end{lem}

\begin{lem}\label{lem:3.2} 

 (5) If $\left\{x_i\right\}_{i=0}^{\infty}$ is an asymptotic average pseudo-orbit of $f$ then there exists a set $J \subset \mathbb{N}$ of density zero such that
$$
\lim _{i \neq J} d\left(x_i, A\right)=0.
$$
\end{lem}

\begin{lem}\label{lem:3.3}

(5) Let $\left\{a_i\right\}_{i=0}^{\infty}$ be a bounded sequence of non-negative real numbers. The following conditions are equivalent:

(1) $\lim \limits_{n \rightarrow \infty} \frac{1}{n} \sum\limits_{i=0}^{n-1} a_i=0$;

(2) There exists a set $J \subset \mathbb{N}$ of density zero such that $\lim \limits_{n \neq J} a_n=0$.
\end{lem}

\textbf{Proof of  Theorem 3.3. }

 For any $i \in \mathbb{N}$, by the compactness of $X$ and continuity of $\mathcal{F}$, let $\left\{z_0, \ldots, z_{2^{i+1}-1}\right\}$ be $\delta_1$-pseudo of $\mathcal{F}\left(\delta_1>0\right)$. Take $u \in X$ and $d\left(z_0, u\right)<\delta_1$, then for $0 \leq j<2^{i+1}$, we get
$$
d\left(F_j(u), z_j\right)<\frac{1}{2^i} .
$$

Let $\left\{x_i\right\}_{i=0}^{\infty} \subset X$ be asymptotic average pseudo-orbit of $\mathcal{F}$ and $Q^{(1)}$ be a density of 0 which satisfies lemma 3.2. Then for sequence $\left\{d\left(f_i\left(x_i\right), x_{i+1}\right)\right\}_{i=0}^{\infty}$, by lemma 3.3, we can find a set $Q^{(2)}$ whose density is 0 such that
$$
\lim _{i \notin Q^{(2)}} d\left(f_i\left(x_i\right), x_{i+1}\right)=0.
$$

Let $J=Q^{(1)} \cup Q^{(2)}$, obviously, the density of $J$ is also 0 . For $n \in \mathbb{N}$, let
$$
J_n=\underset{l \in \Omega_n}{\cup}\left\{l 2^n, \ldots,(l+1) 2^n-1\right\}
$$

\noindent and $J_n$ contain at least one element of $J$, that is
$$
\left\{l 2^n, \ldots,(l+1) 2^n-1\right\} \cap J \neq \varnothing .
$$

\noindent According the way of $J_n$, we get
$$
\begin{gathered}
\lim _{k \notin J_n} d\left(x_k, A\right)=0, \\
\lim _{k \notin J_n} d\left(f_k\left(x_k\right), x_{k+1}\right)=0,
\end{gathered}
$$

\noindent and

$$
\lim _{k \rightarrow+\infty} \sup \frac{1}{k+1} \#\left(J_n \cap[0, k]\right) \leq \lim _{k \rightarrow+\infty} \sup \frac{2^n}{k+1} \#\left(J_n \cap[0, k]\right)=0.
$$

\noindent So, density of every $J_n$ is 0 .

Since
$$
\begin{aligned}
& \lim _{k \neq J} d\left(f_k\left(x_k\right), x_{k+1}\right)=0, \\
& \lim _{k \neq J} d\left(x_k, A\right)=0(i \in \mathbb{N}),
\end{aligned}
$$

\noindent there exists $K_i \in \mathbb{N}$ such that
$$
d\left(f_k\left(x_k\right), x_{k+1}\right)<\delta_i, d\left(x_k, A\right)<\delta_i, \quad \forall k \notin J, k \geq K_i .
$$

Now, we define
$$
R_i=\left\{l 2^{i+1}-1: l \in \Omega_{i+1}, l \geq K\right\}, k \geq K .
$$

\noindent Use lemma 3.1, we can find $J^{\prime} \subseteq \mathbb{N}$. By the way of $R_i$, we get

(i) $\forall i \in \mathbb{N}$, for some $l \in \Omega_{i+1}, m_i=l 2^{i+1}-1$;

(ii) assume $k \notin J^{\prime}, k \geq m_j$, then $d\left(f_k\left(x_k\right), x_{k+1}\right)<\delta_i$;

(iii) assume $k \notin J^{\prime}, k \geq m_j$, then $d\left(x_k, A\right)<\delta_i$.

By the above properties the set $\mathbb{N} / J^{\prime}$ can be written as a union of pairwise disjoint finite sequences of consecutive numbers of non-decreasing length, where the length of each sequence equals $2^i$ for some $i \in \mathbb{N}$. Denote these sequences by $\left\{a_i, \ldots, b_i\right\}$ where
$$
a_0<b_0<a_1<b_1<\ldots
$$

\noindent and put $B=\left\{b_i: i \in \mathbb{N}\right\}$. By the way of $J_n$, we get that if $m_i \leq a_k \leq m_{i+1}$, then

$$
b_k-a_k=2^{i+1} \text {. }
$$

\noindent Note that for every $b_k>m_{i+1}$, we have $b_{k+1}>b_k \geq 2^{i+1}$. So the density of $B$ is 0 .

Take $p \in A$, we  define a sequence $\left\{y_i\right\}_{i=0}^{\infty} \subset A$ be the asymptotic average pseudo-orbit of $\mathcal{F}$. For every $i \in J^{\prime}$, let $y_i=p$, and for every $a_i$, take $y_{a_i}$ such that
$$
d\left(x_{a_i}, y_{a_i}\right)=d\left(x_{a_i}, A\right) .
$$

Let
$$
\left\{y_{a_{i+1}}, \ldots, y_{b_i}\right\}=\left\{F_1\left(y_{a_i}\right), \ldots, F_{b_i-a_i}\left(y_{a_i}\right)\right\} .
$$

\noindent Since $\left\{d\left(f_i\left(y_i\right), y_{i+1}\right)\right\}_{i=0}^{\infty}$ satisfy lemma 3.3 (ii) on $J^{\prime} \cup B$ whose density is 0 , it is an asymptotic average pseudo-orbit of $\left.\mathcal{F}\right|_A$. If $a_i>m_k$, then $\left\{x_j\right\}_{j=a_i}^{b_i}$ is $\delta_k$-pseudo. Hence,
$$
d\left(x_{a_{i+j}}, F_j\left(y_{a_i}\right)\right)<2^{-k}, j=0, \ldots, b_i-a_i .
$$

\noindent Particularly, $\lim _{i \notin J} d\left(x_i, y_i\right)=0$.

Finally, according the asymptotic average shadowing property of $\left.\mathcal{F}\right|_A$, we can select a point $y \in A$ asymptotically shadows $\left\{y_i\right\}_{i=0}^{\infty}$ on average, that is,for some set $C$ of density zero,
$$
\lim _{i \notin C} d\left(F_i(y), y_i\right)=0 .
$$

\noindent Observe the density of $J^{\prime} \cup C$ is 0 , and
$$
\lim _{i \neq J^{\prime} \cup C} \sup d\left(F_i(y), x_i\right) \leq \lim _{i \notin J} \sup d\left(F_i(y), y_i\right)+\lim _{i \neq J^{\prime} \cup C} \sup d\left(y_i, x_i\right)=0.
$$

So by lemma 3.3, $y$ asymptotically shadows $\left\{x_i\right\}_{i=0}^{\infty}$ on average. Hence, $\mathcal{F}$ has asymptotic average shadowing property.

\section{Shadowing For Product Systems}

Let $\left(X, d_1\right)$ and $\left(Y, d_2\right)$ be metric spaces, $\mathcal{F}=\left\{f_n\right\}_{n \in \mathbb{N}}$ and $\mathcal{G}=\left\{g_n\right\}_{n \in \mathbb{N}}$ be time-varying maps on $X$ and $Y$. Define $d$ on $X \times Y$,
$$
d\left\{\left(x_1, y_1\right),\left(x_2, y_2\right)\right\}=\max \left\{d_1\left(x_1, y_1\right), d_2\left(x_2, y_2\right)\right\} .
$$

\noindent Then $d$ is a metric on $X \times Y$.
Let $f: X \rightarrow X, g: Y \rightarrow Y, f$ and $g$ are continuous. Define map,
$$
(f \times g)(x, y)=(f(x), g(y)),(x, y) \in(X, Y),
$$

\noindent then we call $f \times g$ be the product map of $f$ and $g$.

$\operatorname{Let}(F \times G)_K=\left(f_k \times g_k\right)\left(f_{k-1} \times g_{k-1}\right) \ldots\left(f_0 \times g_0\right)$, we get
$$
(\mathcal{F} \times \mathcal{G})_k\left(x_1, x_2\right)=\left(F_k\left(x_1\right), G_k\left(x_2\right)\right) .
$$

\noindent Then we call $\mathcal{F} \times \mathcal{G}=\left\{f_k \times g_k\right\}_{k \in \mathbb{Z}}$ be the family of time- varying mappings of. Hence, $(X \times Y, \mathcal{F} \times \mathcal{G})$ is non-autonomous dynamical systems.

\begin{defn} 
  (2) The sequence $\mathcal{F}=\left\{f_n\right\}_{n \in \mathbb{Z}_{+}}$has $h$-shadowing property if for every $\varepsilon>0$, there exists $\delta>0$ such that for every finite $\delta$-pseudo-orbit $\left\{x_0, x_1, \ldots x_n\right\} \subseteq X$ there is $y \in X$ with $d\left(y, x_0\right)<\varepsilon$ and for all $1 \leq i<n$,
$$
d\left(F_i(y), x_i\right)<\varepsilon, F_n(y)=x_n .
$$

In general, non-autonomous dynamical systems with $h$-shadowing property in compact metric Spaces do not have tracking properties. For example, any double infinite displacement of a finite type with at least two elements.

\end{defn}

\begin{defn} 

(2) Let $\mathcal{F}=\left\{f_n\right\}_{n \in \mathbb{Z}_{+}}$be a time varying map on a metric space $(X, d)$ and $Y$ be a subset of $X$. We say that $\mathcal{F}$ has $s$-limit shadowing property on $Y$, if for every $\varepsilon>0$, there exists $\delta>0$ such that

(i) for every $\delta$-pseudo orbit $\left\{x_n\right\}_{n \in \mathbb{Z}_{+}}$in $Y$, there exists $y \in X \mathcal{E}$-shadows $\left\{x_n\right\}_{n \in \mathbb{Z}_{+}}$;

(ii) $\left\{x_n\right\}_{n \in \mathbb{Z}_{+}}$is a limit pseudo orbit in $\mathrm{Y}$, there exists $y \in X$ limit shadows $\left\{x_n\right\}_{n \in \mathbb{Z}_{+}}$.

In general, systems with limit shadowing properties do not necessarily have $s$-limit shadowing properties.
\end{defn}

\begin{thm}\label{thm:D}
 
 Let $(X, \mathcal{F})$ and $(Y, \mathcal{G})$ be non-autonomous dynamical systems, $d_1$ and $d_2$ be metric respectively. Then $\mathcal{F}$ and $\mathcal{G}$ has $h$-shadowing property if and only if $\mathcal{F} \times \mathcal{G}$ has $h$ shadowing property.
 
\end{thm}

\textbf{Proof of  Theorem 4.1.} 

We first prove the necessity. If $\mathcal{F}$ and $\mathcal{G}$ has $h$-shadowing property, then for any $\varepsilon>0$, there exists $\delta_1>0$ and $\delta_2>0$, such that for finite $\delta_1$-pseudo $\left\{x_0, x_1, \ldots x_n\right\} \subseteq X$ and $\delta_2$-pseudo $\left\{y_0, y_1, \ldots y_n\right\} \subseteq Y$, there are $x \in X$ and $y \in Y$, satisfying $d\left(x, x_0\right)<\varepsilon$ and $d\left(y, y_0\right)<\varepsilon$.

\noindent For all $1 \leq i<n$, we get
$$
\begin{gathered}
d\left(F_i(x), x_i\right)<\varepsilon, F_n(x)=x_n , \\
d\left(G_i(y), y_i\right)<\varepsilon, G_n(y)=y_n .
\end{gathered}
$$

Take $\delta_0=\min \left\{\delta_1, \delta_2\right\},\left\{\left(x_i, y_i\right)\right\}_{i=0}^n$ be the finite $\delta_0$-pseudo orbit of $\mathcal{F} \times \mathcal{G}$. Then for $\forall 1 \leq i<n$,
$$
\begin{aligned}
& d\left((F \times G)_i(x, y),\left(x_i, y_i\right)\right) \\
& \leq d_1\left(F_i(x), x_i\right)+d_2\left(G_i(y), y_i\right) \\
& <2 \varepsilon .
\end{aligned}
$$
and $(F \times G)_n(x, y)=\left(x_n, y_n\right)$. By definition, $\mathcal{F} \times \mathcal{G}$ has $h$-shadowing property.

Now we prove the sufficiency. Let $\left\{\left(x_i, y_i\right)\right\}_{i=0}^n \subseteq X \times Y$ be a finite $\delta$-pseudo orbit of $\mathcal{F} \times \mathcal{G}$. According the $h$-shadowing property of $\mathcal{F} \times \mathcal{G}$, for any $\varepsilon>0$, there exists $\delta>0,(x, y) \in X \times Y$ satisfying
$$
d\left((x, y),\left(x, x_0\right)\right)<\varepsilon, d\left(y, y_0\right)<\varepsilon,
$$

\noindent And for all $1 \leq i<n$, we have
$$
d\left((F \times G)_i(x, y),\left(x_i, y_i\right)\right)<\varepsilon,(F \times G)_n(x, y)=\left(x_n, y_n\right) .
$$

\noindent By the definition of product metric space,
$$
\max \left\{d_1\left(F_i(x), x_i\right), d_2\left(G_i(y), y_i\right)\right\}=d\left((F \times G)_i(x, y),\left(x_i, y_i\right)\right)<\varepsilon .
$$

\noindent So
$$
d_1\left(F_i(x), x_i\right)<\varepsilon, d_2\left(G_i(y), y_i\right)<\varepsilon, F_n(x)=x_n, G_n(y)=y_n .
$$

\noindent Hence, $\mathcal{F}$ and $\mathcal{G}$ has $h$-shadowing property.

\begin{thm}\label{thm:E}
 
Let $(X, \mathcal{F})$ and $(Y, \mathcal{G})$ be non-autonomous dynamical systems, $d_1$ and $d_2$ be metric respectively. Then $\mathcal{F}$ and $\mathcal{G}$ has $s$-limit shadowing property if and only if $\mathcal{F} \times \mathcal{G}$ has $s$-limit shadowing property.

\end{thm}

\textbf{Proof of Theorem 4.2.}

We first prove the necessity. By the $s$-limit shadowing property of $\mathcal{F}$, for any $\varepsilon>0$, there exists $\delta_1>0$,such that for every $\delta_1$-pseudo orbit $\left\{x_n\right\}_{n \in \mathbb{Z}_{+}}$in $X$, there is $x \in X$ such that for all $n \geq 0$,
$$
d\left(F_n(x), x_n\right)<\varepsilon .
$$

\noindent And if $\left\{x_n\right\}_{n \in \mathbb{Z}}$ is the limit pseudo orbit, then $d\left(F_n(x), x_n\right) \rightarrow 0(n \rightarrow \infty)$. Similarly, by the $s$-limit shadowing property of $\mathcal{G}$, there exists $\delta_2>0$ such that for every $\delta_2$-pseudo orbit $\left\{y_n\right\}_{n \in \mathbb{Z}}$ in $Y$, there is $y \in Y$ such that for all $n \geq 0$,

$$
d\left(G_n(y), y_n\right)<\varepsilon .
$$

If $\left\{y_n\right\}_{n \in \mathbb{Z}_{+}}$is the limit pseudo orbit, then $d\left(G_n(y), y_n\right) \rightarrow 0(n \rightarrow \infty)$.

Take $\delta_0=\min \left\{\delta, \delta^{\prime}\right\},\left\{\left(x_n, y_n\right)\right\}_{n \geq 0}$ be $\delta_0$-pseudo orbit of $\mathcal{F} \times \mathcal{G}$, then there is $(x, y) \in X \times Y$ such that for all $n \geq 0$, we get
$$
\begin{aligned}
& d\left((F \times G)_n(x, y),\left(x_n, y_n\right)\right) \\
& \leq d_1\left(F_n(x), x_n\right)+d_2\left(G_n(y), y_n\right) \\
& <2 \varepsilon .
\end{aligned}
$$

\noindent And if $\left\{\left(x_n, y_n\right)\right\}_{n \in \mathbb{Z}}$ is the limit pseudo orbit in $X \times Y$. Obviously, we have
$$
d\left((F \times G)_n(x, y),\left(x_n, y_n\right)\right) \rightarrow 0(n \rightarrow \infty) .
$$

\noindent So $\mathcal{F} \times \mathcal{G}$ has $s$-limit shadowing property.

Now we prove the sufficiency. Let $\left\{\left(x_n, y_n\right)\right\}_{n \geq 0} \subseteq X \times Y$ be a $\delta$-pseudo orbit of $\mathcal{F} \times \mathcal{G}$. By the $s$-limit shadowing property of $\mathcal{F} \times \mathcal{G}$, there exists $(x, y) \in X \times Y$ such that for all $n \geq 0$,
$$
d\left((F \times G)_n(x, y),\left(x_n, y_n\right)\right)<\varepsilon .
$$

\noindent Hence,

$$
\max \left\{d_1\left(F_n(x), x_n\right), d_2\left(G_n(y), y_n\right)\right\}=d\left((F \times G)_n(x, y),\left(x_n, y_n\right)\right)<\varepsilon .
$$

If $\left\{\left(x_i, y_i\right)\right\}_{i \geq 0}$ is limit pseudo orbit, that is,
$$
\lim _{n \rightarrow+\infty} d\left((f \times g)_{n+1}\left(x_n, y_n\right),\left(x_{n+1}, y_{n+1}\right)\right)=0,
$$

\noindent then

$$
\lim _{n \rightarrow+\infty} d\left((F \times G)_n(x, y),\left(x_n, y_n\right)\right)=0 .
$$

\noindent So for any $\varepsilon>0$, there exists $\mathrm{N}>0$, if $n \geq 0$,
$$
d\left((f \times g)_{n+1}\left(x_n, y_n\right),\left(x_{n+1}, y_{n+1}\right)\right)<\varepsilon .
$$

\noindent Hence,
$$
d_1\left(f_{n+1}\left(x_n\right), x_{n+1}\right)<\varepsilon, d_2\left(g_{n+1}\left(y_n\right), y_{n+1}\right)<\varepsilon \text {. }
$$

\noindent Therefore,
$$
\lim _{n \rightarrow+\infty} d_1\left(f_{n+1}\left(x_n\right), x_{n+1}\right)=0,
$$

$$
\lim _{n \rightarrow+\infty} d_2\left(g_{n+1}\left(y_n\right), y_{n+1}\right)=0 .
$$
$\left\{x_n\right\}_{n \geq 0}$ and $\left\{y_n\right\}_{n \geq 0}$ is limit pseudo orbit of $\mathcal{F}$ and $\mathcal{G}$ respectively. Similarly, we also get
$$
\begin{aligned}
& \lim _{n \rightarrow+\infty} d_1\left(F_n(x), x_n\right)=0, \\
& \lim _{n \rightarrow+\infty} d_2\left(G_n(y), y_n\right)=0 .
\end{aligned}
$$

\noindent By definition, $\mathcal{F}$ and $\mathcal{G}$ has $s$-limit shadowing property.

The above theorems also hold for asymptotic average shadowing, limit shadowing, avergage shadowing, period shadowing and Lipschitz shadowing.\\


\bibliography{REFFERENCE}\end{document}